\input amstex\documentstyle{amsppt}  
\pagewidth{12.5cm}\pageheight{19cm}\magnification\magstep1
\topmatter
\title Fourier transform as a triangular matrix\endtitle
\author G. Lusztig\endauthor
\address{Department of Mathematics, M.I.T., Cambridge, MA 02139}\endaddress
\thanks{Supported by NSF grant DMS-1855773.}\endthanks
\endtopmatter   
\document

\define\si{\sim}

\define\sqc{\sqcup}

\define\lb{\linebreak}

\define\bin{\binom}
\define\op{\oplus}
   
\define\part{\partial}
\define\emp{\emptyset}

\define\m{\mapsto}
\define\do{\dots}

\define\sub{\subset}    

\define\T{\times}
\define\ti{\tilde}
\define\nl{\newline}

\define\tr{\text{\rm tr}}

\redefine\b{\beta}

\define\g{\gamma}
\redefine\d{\delta}
\define\e{\epsilon}

\define\ps{\psi}
\define\r{\rho}
\define\s{\sigma}
\redefine\t{\tau}
\define\th{\theta}
\define\k{\kappa}
\redefine\l{\lambda}
\define\z{\zeta}

\redefine\G{\Gamma}

\redefine\L{\Lambda}
\define\Ph{\Phi}

\redefine\ss{\bold s}

\define\BB{\bold B}
\define\CC{\bold C}

\define\FF{\bold F}

\define\QQ{\bold Q}
\define\RR{\bold R}

\define\ZZ{\bold Z}

\define\cf{\Cal F}
\define\cg{\Cal G}
\define\ch{\Cal H}
\define\ci{\Cal I}

\head Introduction\endhead
\subhead 0.1\endsubhead
Let $V$ be a vector space of finite even dimension $D=2d\ge0$ over the field $\FF_2$ with $2$ elements 
with a fixed nondegenerate symplectic form $(,):V\T V@>>>\FF_2$. Let $[V]$ be the $\CC$-vector space 
of functions $V@>>>\CC$ and let $[V]_\ZZ$ be the subgroup of $[V]$ consisting of the functions 
$V@>>>\ZZ$. For $f\in[V]$ the Fourier transform $\Ph(f)\in[V]$ is defined by 
$\Ph(f)(x)=2^{-d}\sum_{y\in V}(-1)^{(x,y)}f(y)$. Now 
$\Ph:[V]@>>>[V]$ is a linear involution whose trace is $2^{-d}\sum_{x\in V}1=2^d$. Hence $\Ph$ 
has $2^{D-1}+2^{d-1}$ eigenvalues equal to $1$ and $2^{D-1}-2^{d-1}$ eigenvalues equal to $-1$. 
Here is one of our main results.

\proclaim{Theorem 0.2} There exists a $\ZZ$-basis $\b$ of $[V]_\ZZ$ consisting of
characteristic functions of certain explicit isotropic subspaces of $V$ such that the matrix of 
$\Ph:[V]@>>>[V]$ with respect to $\b$ is upper triangular (with diagonal entries $\pm1$) for a 
suitable order on $\b$.
\endproclaim
Assume for example that $D=2$. For $x\in V$ let $f_x\in[V]$ be the function whose value at $y\in V$ is 
$1$ if $y=x$ and $0$ if $y\ne x$.
Let $\b$ be the $\ZZ$-basis of $V_\ZZ$ consisting of $f'_0=f_0$ and of $f'_x=f_0+f_x$ for $x\in V-\{0\}$. 
We have $\Ph(f'_0)=-f'_0+(1/2)\sum_{x\in V-\{0\}}f'_x$ and $\Ph(f'_x)=f'_x$ for $x\in V-\{0\}$. Thus,
the matrix of $\Ph:[V]@>>>[V]$ with respect to $\b$ is upper triangular (with diagonal entries 
$-1,1,1,1$). 

The proof of the theorem is given in \S1; we take $\b$ to be the new basis $\cf(V)$ of $[V]$ defined 
in \cite{L19}. In \S2 we compute explicitly the signs $\pm1$ appearing in the theorem for this $\b$.
In \S3 we give some tables for $\b=\cf(V)$. In \S4 we show that $\cf(V)$ has a certain dihedral 
symmetry which was not apparent in \cite{L19}. In \S5 we show that the theorem is a special case of 
a result which applies to any two-sided cell in an irreducible Weyl group.

\subhead 0.3\endsubhead
{\it Notation.} For $a,b$ in $\ZZ$ we set $[a,b]=\{z\in\ZZ;a\le z\le b\}$. 
For a finite set $Y$ let $|Y|$ be the cardinal of $Y$. 

\head 1. Proof of Theorem 0.2\endhead
\subhead 1.1\endsubhead
When $D\ge2$ we 
fix a subset $\{e_i;i\in[1,D+1]\}\sub V$ such that for $i\ne j$ in $[1,D+1]$ we have $(e_i,e_j)=1$ 
if $i-j=\pm1\mod D+1$, $(e_i,e_j)=0$ if $i-j\ne\pm1\mod D+1$. (Such a subset exists and is unique up to 
the action of some isometry of $(,)$.) We say that this subset is a {\it circular basis} of $V$. We 
must have $e_1+e_2+\do+e_{D+1}=0$ and any $D$ elements of $\{e_i;i\in[1,D+1]\}$ form a basis of $V$. 
For any $I\sub[1,D+1]$ let $e_I=\sum_{i\in I}e_i\in V$. 
When $D\ge2$ (resp. $D\ge4$) we denote by $V'$ (resp. $V''$) an $\FF_2$-vector space with a nondegenerate
symplectic form $(,)$. When $D\ge4$ (resp. $D\ge6$) we assume that $V'$ (resp. $V''$) has a given 
circular basis $\{e'_i;i\in[1,D-1]\}$ (resp. $\{e''_i;i\in[1,D-3]\}$).

When $D\ge2$, for any $i\in[1,D+1]$ there is a unique linear map $\t_i:V'@>>>V$ such that 
$\t_i=0$ for $D=2$, while for $D\ge4$, the sequence $\t_i(e'_1),\t_i(e'_2),\do,\t_i(e'_{D-1})$ is:

$e_1,e_2,\do,e_{i-2},e_{i-1}+e_i+e_{i+1},e_{i+2},e_{i+3},\do,e_D,e_{D+1}$ (if $1<i\le D$),

$e_3,e_4,\do,e_D,e_{D+1}+e_1+e_2$  if $i=1$

$e_2,e_3,\do,e_{D-1},e_D+e_{D+1}+e_1$ if $i=D+1$.
\nl
This map is injective and compatible with $(,)$. Similarly, when $D\ge4$, for any $i\in[1,D-1]$ there is a 
unique linear map $\t'_i:V''@>>>V'$ such that $\t'_i=0$ for $D=4$, while for $D\ge6$, the sequence 
$\t'_i(e''_1),\t'_i(e''_2),\do,\t'_i(e''_{D-3})$ is:

$e'_1,e'_2,\do,e'_{i-2},e'_{i-1}+e'_i+e'_{i+1},e'_{i+2},e'_{i+3},\do,e'_{D-2},e'_{D-1}$ (if $1<i\le D-2$),

$e'_3,e'_4,\do,e'_{D-2},e'_{D-1}+e'_1+e'_2$  if $i=1$

$e'_2,e'_3,\do,e'_{D-3},e'_{D-2}+e'_{D-1}+e'_1$ if $i=D-1$.
\nl
This map is injective and compatible with $(,)$. Note that 

(a) if $D\ge2$, then $\t_i(V')$ is a complement of the line $\FF_2e_i$ in $\{x\in V;(x,e_i)=0\}$.
\nl
Assuming that $D\ge4$ and $i\in[1,D-2]$, we show:

(b) $\t_{D+1}\t'_i=\t_j\t'_{D-1}$ where $j=i+1$ if $1<i\le D-2$ and $j=i$ if $i=1$.
\nl
If $D=4$ the result is trivial. Assume now that $D\ge6$.
Assume first that $1<i\le D-2$.  Both sequences 
$$(\t_{D+1}\t'_i(e''_1),\t_{D+1}\t'_i(e''_2),\do,\t_{D+1}\t'_i(e''_{D-3}))$$
$$(\t_{i+1}\t'_{D-1}(e''_1),\t_{i+1}\t'_{D-1}(e''_2),\do,\t_{i+1}\t'_{D-1}(e''_{D-3}))$$
are equal to
$$(e_2,e_3,\do,e_{i-1},e_i+e_{i+1}+e_{i+2},e_{i+3},e_{i+4},\do,e_{D-1},e_D+e_{D+1}+e_1)$$
if $1<i<D-2$ and to
$$(e_2,e_3,\do,e_{D-3},e_{D-2}+e_{D-1}+e_D+e_{D+1}+e_1)$$
if $i=D-2$.
Next we assume that $i=1$.  Both sequences 
$$(\t_{D+1}\t'_i(e''_1),\t_{D+1}\t'_i(e''_2),\do,\t_{D+1}\t'_i(e''_{D-3}))$$
$$(\t_i\t'_{D-1}(e''_1),\t_i\t'_{D-1}(e''_2),\do,\t_i\t'_{D-1}(e''_{D-3}))$$
are equal to
$$(e_4,e_5,\do,e_{D-2},e_D+e_{D+1}+e_1+e_2+e_3).$$
This proves (b).

In the setup of (b) we show that for a subspace $E''\sub V''$ we have
$$\t_{D+1}(\t'_i(E'')\op\FF_2e'_i)\op\FF_2e_{D+1}=\t_j(\t'_{D-1}(E'')\op\FF_2e'_{D-1})\op\FF_2e_j.
\tag c$$
Using (b) it is enough to show that
$$\FF_2\t_{D+1}(e'_i)\op\FF_2e_{D+1}=\FF_2\t_j(e'_{D-1})\op\FF_2e_j$$
or that
$$\FF_2e_{i+1}\op\FF_2e_{D+1}=\FF_2e_{D+1}\op\FF_2e_{i+1}$$
if $i>1$ and
$$\FF_2e_2\op\FF_2e_{D+1}=\FF_2(e_{D+1}+e_1+e_2)+\FF_2e_1$$
if $i=1$. This is clear.

\subhead 1.2\endsubhead
If $D\ge2$, for any $k\in[0,d]$ let $E_k$ be the subspace of $V$ with basis

$\{e_{[1,D]},e_{[2,D-1]},\do,e_{[k,D+1-k]}\}$.
\nl
When $D=0$ we set $E_0=0\sub V$. If $D\ge4$ and $k\in[0,d-1]$ let $E'_k$ be the subspace of $V'$ with basis

$\{e'_{[1,D-2]},e'_{[2,D-3]},\do,e'_{[k,D-1-k]}\}$
\nl
where for any $I'\sub[1,D-1]$ we set $e'_{I'}=\sum_{i\in I'}e'_i\in V'$. When $D=2$ we set $E'=0\sub V'$.

Following  \cite{L19}, we define a collection $\cf(V)$ of subspaces of $V$ by induction on $D$.
If $D=0$, $\cf(V)$ consists of the subspace $\{0\}$. If $D\ge2$, a subspace $E$ of $V$ is in $\cf(V)$
if either

(i) there exists $i\in[1,D]$ and $E'\in\cf(V')$ such that $E=\t_i(E')\op\FF_2 e_i$, or

(ii) there exists $k\in[0,d]$ such that $E=E_k$.
\nl
We now define a collection $\cf'(V)$ of subspaces of $V$ by induction on $D$.
If $D=0$, $\cf'(V)$ consists of the subspace $\{0\}$. If $D\ge2$, a subspace $E$ of $V$ is in 
$\cf'(V)$ if either $E=0$ or if

(iii) there exists $i\in[1,D+1]$ and $E'\in\cf(V')$ such that $E=\t_i(E')\op\FF_2 e_i$.

\proclaim{Lemma 1.3} We have $\cf(V)=\cf'(V)$.
\endproclaim
We argue by induction on $D$. If $D=0$ the result is obvious. Assume that $D\ge2$. We show that 

(a) $\cf'(V)\sub\cf(V)$. 
\nl
Let $E\in\cf'(V)$. If $E=0$ then clearly $E\in\cf(V)$. Thus we can 
assume that $E=\t_i(E')\op\FF_2 e_i$ for some $i\in[1,D+1]$ and some $E'\in\cf'(V)$. By the
induction hypothesis we have $E'\in\cf(V)$. If $i\in[1,D]$ then by definition we have $E\in\cf(V)$.
Thus we can assume that $i=D+1$. If $E'=0$ then $E=\FF_2e_{D+1}=\FF_2e_{[1,D]}=E_1\in\cf(V)$.
Thus we can assume that $E'\ne0$ so that $D\ge4$. Since $E'\in\cf(V')$ we have
$E'=\t'_h(E'')\op\FF_2 e'_h$ for some $h\in[1,D-2]$ and some $E''\in\cf(V'')$. Thus we have
$$E=\t_{D+1}(\t'_h(E'')\op\FF_2e'_h)\op\FF_2e_{D+1}
=\t_{h'}(E_1)\op\FF_2e_{h'}$$
where $E_1=\t'_{D-1}(E'')\op\FF_2e_{D-1}$
(we have used 1.1(c)); here $h'=h+1$ if $h>1$ and $h'=h$ if $h=1$.
By the definition of $\cf'(V')$ we have $E_1\in\cf'(V')$ hence
$E_1\in\cf(V')$, by the induction hypothesis.
It follows that $\t_{h'}(E_1)\op\FF_2e_{h'}\in\cf(V)$, so that $E\in\cf(V)$. This proves (a).

We show that 

(b) $\cf(V)\sub\cf'(V)$. 
\nl
Let $E\in\cf(V)$. Assume first that $E=E_k$ for some $k\in[1,d]$.
From the definition we have
$E_k=\t_{D+1}(E'_{k-1})\op\FF_2e_{D+1}$.
We have $E'_{k-1}\in\cf(V')$ hence by the induction hypothesis we have
$E'_{k-1}\in\cf'(V')$ and using the definition we have $E_k\in\cf'(V)$.
If $E=E_0$ then $E=0$ so that again $E\in\cf(V)$.
Next we assume that $E$ is not of the form $E_k$ with $k\in[0,d]$.
We can find $i\in[1,D]$ and $E'\in\cf(V')$ such that $E=\t_i(E')\op\FF_2e_i$.
By the induction hypothesis we have $E'\in\cf'(V')$. From the definition we have $E\in\cf'(V)$.
This proves (b).

\subhead 1.4\endsubhead
For any subset $X\sub V$ let $\ps_X\in[V]$ be the function such that $\ps_X(x)=1$ if $x\in X$,
$\ps_X(x)=0$ if $x\in V-X$. According to \cite{L19},

(a) {\it $\{\ps_E;E\in\cf(V)\}$ is a $\ZZ$-basis of $[V]_\ZZ$.}
\nl
Using Lemma 1.3, we deduce:

(b) {\it $\{\ps_E;E\in\cf'(V)\}$ is a $\ZZ$-basis of $[V]_\ZZ$.}
\nl
We will no longer distinguish between $\cf(V)$ and $\cf'(V)$.

\subhead 1.5\endsubhead
Assume that $D\ge2$. Let $[V']$, $\Ph':[V']@>>>[V']$ be the analogues of $[V],\Ph:[V]@>>>[V]$ when
$V$ is replaced by $V'$. For $X'\sub V'$ let $\ps'_{X'}\in[V']$ be the function such that 
$\ps'_{X'}(y)=1$ if $y\in X'$, $\ps'_{X'}(x)=0$ if $y\in V'-X'$. 

For $i\in[1,D+1]$ there is a unique linear map $z_i:[V']@>>>[V]$ such that 
$z_i(\ps'_y)=\ps_{\t_i(y)}+\ps_{\t_i(y)+e_i}$ for all $y\in V'$.
If $E'$ is a subspace of $V'$ we have $z_i(\ps'_{E'})=\ps_{\t_i(E')\op\FF_2e_i}$. We show:

(a) For $f\in[V']$ we have $\Ph(z_i(f))=z_i(\Ph'(f))$.
\nl
We can assume that $f=\ps'_y$ with $y\in V'$. We have 
$$\align&z_i(\Ph'(f))=2^{-d+1}\sum_{y_1\in V'}(-1)^{(y,y_1)}z_i(\ps'_{y_1})\\&=
2^{-d+1}\sum_{y_1\in V'}(-1)^{(y,y_1)}(\ps_{\t_i(y_1)}+\ps_{\t_i(y_1)+e_i}),\endalign$$
$$\align&\Ph(z_i(f))=\Ph(\ps_{\t_i(y)}+\ps_{\t_i(y)+e_i})=
2^{-d}\sum_{x\in V}((-1)^{(\t_i(y),x)}+(-1)^{(\t_i(y)+e_i,x)})\ps_x\\&=
2^{-d+1}\sum_{x\in V;(e_i,x)=0}(-1)^{(\t_i(y),x)}\ps_x.\endalign$$
In the last sum $x$ can be written uniquely as $x=\t_i(y_1)+ce_i$ with $y_1\in V',c\in\FF_2$.
Thus
$$\Ph(z_i(f))=2^{-d+1}\sum_{y_1\in V',c\in\FF_2}(-1)^{(\t_i(y),\t_i(y_1)+ce_1)}\ps_{\t_i(y_1)+ce_i}$$
which is equal to $z_i(\Ph'(f))$. This proves (a).

For $E\in\cf(V)$ we write 

(b) $\Ph(\ps_E)=\sum_{E_1\in\cf(V)}c_{E,E_1}\ps_{E_1}$ 
\nl
with $c_{E,E_1}\in\CC$ are uniquely determined. (We use 1.4(b).)

\proclaim{Lemma 1.6} Let $E\in\cf(V),E_1\in\cf(V)$ be such that $c_{E,E_1}\ne0$.
Then either $E_1=E$ or $|E_1|>|E|$.
\endproclaim
We argue by induction on $D$. If $D=0$ the result is obvious. Assume now that $D\ge2$.
If $E=0$, the result is obvious since for any $E_1\in\cf(V)$ we have either $E_1=E$ or $|E_1|>|E|$.
Assume now that $E\ne0$.
We can find $i\in[1,D+1]$ and $E'\in\cf(V')$ such that $E=\t_i(E')\op\FF_2e_{\ti i}$.
Recall from 1.5 that $z_i(\ps'_{E'})=\ps_{\t_i(E')\op\FF_2e_i}=\ps_E$.
By the induction hypothesis we have
$$\Ph'(\ps'_{E'})=c'_{E',E'}\ps'_{E'}+\sum_{E'_1\in\cf(V');|E'_1|>|E'|}c'_{E',E'_1}\ps'_{E'_1}$$
 with $c'_{E',E'}\in\CC$, $c'_{E',E'_1}\in\CC$. Applying $z_i$ and using 1.5(a) we deduce
$$\align&\Ph(z_i(\ps'_{E'}))=c_{E',E'}z_i(\ps'_{E'})+\sum_{E'_1\in\cf(V');|E'_1|>|E'|}c'_{E',E'_1}
z_i(\ps'_{E'_1})\\&=
c_{E',E'}\ps_E\sum_{E'_1\in\cf(V');|E'_1|>|E'|}c'_{E',E'_1}\ps_{\t_i(E'_1)\op\FF_2 e_i}\endalign$$
and the result follows in this case since for $E'_1$ in the last sum we have
$$|\t_i(E'_1)\op\FF_2 e_i|=|E'_1|+1>|E'|+1=|E|.$$
This completes the proof of the lemma.

\subhead 1.7\endsubhead
We prove Theorem 0.2. By results of \cite{L19}, the basis 1.4(b) of $[V]$ is a $\ZZ$-basis of 
$[V]_\ZZ$. By 1.6, the matrix of $\Ph$ with respect to the basis 1.4(b) is upper triangular for a 
suitable order on the basis. The diagonal entries of this matrix are necessarily $\pm1$ since 
$\Ph^2=1$. This completes the proof.

\head 2. Sign computation\endhead
\subhead 2.1\endsubhead
In this subsection we assume that $D\ge2$.
Let $E\in\cf(V)$. According to \cite{L19} there is a unique basis
$b_E$ of $E$ which consists of vectors of the form $e_I$ with $I$ of the form $[a,b]$
with $a\le b$ in $[1,D]$. Let $n_E$ be the number of vectors $e_I\in b_E$ such that $|I|$ is even.

For $k\in[0,d]$ let $\cf_k(V)$ (resp. $\cf^k(V)$) be the set of all $E\in\cf(V)$ such that 
$\dim(E)=k$ (resp. $n_E=k$). The following equality follows from \cite{L19, 1.27(a)}:

(a) {\it $|\cf^k(V)|=\bin{D+1}{d-k}$.}
\nl
If $E\in\cf(V)$ we denote by $E^!$ the subspace of $E$ spanned by the vectors $e_I\in b_E$
such that $|I|$ is odd; we have $E^!\in\cf^0(V)$.
One can show that for any $E_1\in\cf^0(V)$, the set $\{E\in\cf(V);E^!=E_1\}$
consists of $k+1$ subspaces
$E_1(0),E_1(1),\do,E_1(k)$ where $k\in[0,d]$ is defined by $E_1\in\cf_{d-k}(V)$ and 
$E_1(0)\in\cf^{d-k}(V)$, $E_1(1)\in\cf^{d-k+1}(V)$, $\do$, $E_1(k)\in\cf^d(V)$. 
(This is illustrated in the tables in \S3.) We have $E_1(0)=E_1$. 
We define an involution $\k:\cf(V)@>>>\cf(V)$ by the requirement that 
$\k(E_1(0))=E_1(k)$, $\k(E_1(1))=E_1(k-1)$, $\do$, $\k(E_1(k))=E_1(0)$ for any 
$E_1\in\cf^0\cap\cf_{d-k}(V)$. This involution restricts to a bijection
$$\cf^k(V)@>\si>>\cf_{d-k}(V)$$
for any $k\in[0,d]$. Using this bijection and (a) we deduce

(b) {\it $|\cf_k(V)|=\bin{D+1}{k}$ for $k\in[0,d]$.}

\subhead 2.2\endsubhead
For any integer $N$ we set $\d(N)=(-1)^{N(N+1)/2}$. We have the following identity:

(a) $\sum_{k\in[0,d]}\d(d-k)\bin{D+1}{k}=2^d$.
\nl
We prove (a) by induction on $D$. If $D=0$ the result is obvious. Assume now that $D\ge2$. 
We must show that
$$\bin{2d+1}{d}-\bin{2d+1}{d-1}-\bin{2d+1}{d-2}+\bin{2d+1}{d-3}+\bin{2d+1}{d-4}-\bin{2d+1}{d-5}-\do=
2^d$$
or that
$$\align&(\bin{2d}{d}+\bin{2d}{d-1})-(\bin{2d}{d-1}+\bin{2d}{d-2})-(\bin{2d}{d-2}+\bin{2d}{d-3})\\&+
(\bin{2d}{d-3}+\bin{2d}{d-4})+(\bin{2d}{d-4}+\bin{2d}{d-5})-(\bin{2d}{d-5}+\bin{2d}{d-6})-\do\\&
=2^d\endalign$$
or that
$$\bin{2d}{d}-2\bin{2d}{d-2}+2\bin{2d}{d-4}-2\bin{2d}{d-6}+\do=2^d$$
or that 
$$\align&(\bin{2d-1}{d}+\bin{2d-1}{d-1})-2(\bin{2d-1}{d-2}+\bin{2d-1}{d-3}\\&+
2(\bin{2d-1}{d-4}+\bin{2d-1}{d-5})-\do=2^d\endalign$$
or that
$$2\bin{2d-1}{d-1}-2\bin{2d-1}{d-2}-2\bin{2d-1}{d-3}+2\bin{2d-1}{d-4}+2\bin{2d-1}{d-5}-\do=2^d.$$
But this is known from the induction hypothesis. This proves (a).

\subhead 2.3\endsubhead
The following result describes the diagonal entries of the upper triangular matrix in 1.7. 

\proclaim{Proposition 2.4} Let $E\in\cf(V)$ and let $c_{E,E}$ be as in 1.5(b). We have 
$c_{E,E}=\d(d-\dim E)$.
\endproclaim
We argue by induction on $D$. If $D=0$ the result is obvious. Assume now that $D\ge2$. Assume first 
that $E\ne0$. We can find $i\in[1,D+1]$ and $E'\in\cf(V')$ such that $E=\t_i(E')\op\FF_2e_{\ti i}$.
By the proof of 1.6 we have $c_{E,E}=c'_{E',E'}$ (notation of 1.6). The proposition applies
to $c'_{E',E'}$ by the induction hypothesis. The desired result for $E$ follows since
$d-\dim E=d-1-\dim E'$.
We now assume that $E=0$. The trace of $\Ph$ is equal to
$\sum_{E_1\in\cf(V)}c_{E_1,E_1}$ and on the other hand is equal to $2^d$ (see 0.1).
Thus we have $\sum_{E_1\in\cf(V)}c_{E_1,E_1}=2^d$.
In the last sum all terms with $E_1\ne0$ are already known. Hence the term with $E_1=0$ is determined
by the last equality. Thus to prove the proposition it is enough to verify the identity
$$\sum_{E_1\in\cf(V)}\d(d-\dim E_1)=2^d$$
or equivalently
$$\sum_{k\in[0,d]}|\cf_k(V)|\d(d-k)=2^d.$$
This follows from 2.1(b), 2.2(a). This completes the proof.

\head 3. Tables\endhead
\subhead 3.1\endsubhead
In this section we assume that $D\ge2$.
Let $E\in\cf(V)$. Recall that the basis $b_E$ consists of certain vectors $e_I$
where $I$ is of the form $[a,b]$ with $a\le b$ in $[1,D]$. We have $e_I=e_{I'}$
where $I'\sub[1,D+1]$ is defined by $I'=I$ if $|I|$ is odd and $I'=[1,D+1]-I$ if $I$ is even.
Note that $|I'|$ is always odd. Now $E$ is completely described by the list of all subsets $I'$
defined as above. In the following three sections we describe each $E\in\cf(V)$ as a list
of such $I'$ assuming that $D$ is $2,4$ or $6$. 
(This list is more symmetric than the corresponding list of the $I$ which is given in \cite{L19}.)
In each of these tables each horizontal line represents the various $E_1(0),E_1(1),\do,E_1(k)$
with a fixed $E_1\in\cf^0(V)$ as in 2.1. For example the second line $<1>,<1,512>$ in 3.3
represents two subspaces in $\cf(V)$; one spanned by $e_1$ and the other spanned by $e_1$ and
$e_5+e_1+e_2$.

\subhead 3.2. The table for $D=2$\endsubhead

$\emp,<3>$

$<1>$

$<2>$.

\subhead 3.3. The table for $D=4$\endsubhead

$\emp,<5>, <5,451>$

$<1>,<1,512>$

$<2>,<2,5>$

$<3>,<3,5>$

$<4>,<4,345>$

$<1,3>$

$<1,4>$

$<2,4>$ 

$<2,123>$

$<3,234>$.

\subhead 3.3. The table for $D=6$\endsubhead

$\emp,<7>,<7,671>,<7,671,56712>$

$<1>,<1,712>,<1,712,67123>$

$<2>,<2,7>,<2,7,67123>$

$<3>,<3,7>,<3,7,671>$

$<4>,<4,7>,<4,7,671>$

$<5>,<5,7>,<5,7,45671>$

$<6>,<6,567>,<6,567,45671>$

$<1,3>,<1,3,71234>$

$<1,4>,<1,4,712>$

$<1,5>,<1,5,712>$

$<1,6>,<1,6,56712>$

$<2,4>,<2,4,7>$

$<2,5>,<2,5,7>$

$<2,6>,<2,6,567>$

$<3,5>,<3,5,7>$

$<3,6>,<3,6,567>$

$<4,6>,<4,6,34567>$

$<2,123>,<2,123,71234>$

$<3,234>,<3,7,234>$

$<4,345>,<4,7,345>$

$<5,456>,<5,456,34567>$

$<1,3,5>$

$<1,3,6>$

$<1,4,6>$

$<2,4,6>$

$<1,4,345>$

$<1,5,456>$

$<2,5,123>$

$<2,5,456>$

$<2,6,123>$

$<3,6,234>$

$<2,4,12345>$

$<3,5,23456>$

$<3,234,12345>$

$<4,345,23456>$.

\head 4. Dihedral symmetry\endhead
\subhead 4.1\endsubhead
There is a unique linear map $R:V@>>>V$ such that if $D=0$ we have $R=0$ while if $D\ge2$,
$R(e_1),R(e_2),\do,R(e_{D+1})$ is 
$e_2,e_3,\do,e_{D+1},e_1$. If $D\ge2$, there is a unique linear map $R':V'@>>>V'$ such that 
if $D=2$ we have $R'=0$ while if $D\ge4$,
$R'(e'_1),R'(e'_2),\do,R'(e'_{D-1})$ is $e'_2,e'_3,\do,e'_{D-1},e'_1$. From the definitions
we see that if $D\ge2$, $i\in[1,D+1]$ we have

(a) $R\t_i=\t_{i+1}R':V'@>>>V$ if $i\in[1,D]$, $R\t_i=\t_1:V'@>>>V$ if $i=D+1$.

\subhead 4.2\endsubhead
Let $E\in\cf(V)$. We show:

(a) $R(E)\in\cf(V)$.
\nl
We argue by induction on $D$. If $D=0$ the result is obvious. Assume that $D\ge2$. If $E=0$
we have $R(E)=0$ and the result is clear. Assume now that $E\ne0$. We can find $i\in[1,D+1]$ and 
$E'\in\cf(V')$ such that $E=\t_i(E')\op\FF_2e_i$. Applying $R$ we deduce 
$R(E)=R\t_i(E')\op\FF_2e_{i+1}$ if $i\in[1,D]$, $R(E)=R\t_i(E')\op\FF_2e_1$ if $i=D+1$. Using 
4.1(a) we deduce $R(E)=\t_{i+1}R'(E')\op\FF_2e_{i+1}$ if $i\in[1,D]$, $R(E)=\t_1(E')\op\FF_2e_1$ 
if $i=D+1$. By the induction hypothesis we have $R'(E')\in\cf(V')$. 
It follows that $R(E)\in\cf(V)$, as required.

\subhead 4.3\endsubhead
There is a unique linear map $S:V@>>>V$ such that if $D=0$ we have $S=0$, while if $D\ge2$ we have
      $S(e_i)=e_{D+1-i}$ if $i\in[1,D]$,
$S(e_{D+1})=e_{D+1}$. If $D\ge2$, there is a unique linear map $S':V'@>>>V'$ such that
if $D=2$ we have $S'=0$ while if $D\ge4$ we have 
$S'(e_i)=e_{D-1-i}$ if $i\in[1,D-2]$, $S'(e_{D-1})=e_{D-1}$.
From the definitions we see that if $D\ge2$, $i\in[1,D+1]$ we have

(a) $S\t_i=\t_{D+1-i}S':V'@>>>V$ if $i\in[1,D]$, $S\t_i=\t_iS':V'@>>>V$ if $i=D+1$.

\subhead 4.4\endsubhead
Let $E\in\cf(V)$. We show:

(a) $S(E)\in\cf(V)$.
\nl
We argue by induction on $D$. If $D=0$ the result is obvious. Assume that $D\ge2$. If $E=0$
we have $S(E)=0$ and the result is clear. Assume now that $E\ne0$. We can find $i\in[1,D+1]$ and 
$E'\in\cf(V')$ such that $E=\t_i(E')\op\FF_2e_i$. Applying $S$ we deduce 
$S(E)=S\t_i(E')\op\FF_2e_{D+1-i}$ if $i\in[1,D]$, $S(E)=S\t_i(E')\op\FF_2e_i$ if $i=D+1$. Using 
4.3(a) we deduce $S(E)=\t_{D+1-i}S'(E')\op\FF_2e_{D+1-i}$ if $i\in[1,D]$, 
$S(E)=\t_iS'(E')\op\FF_2e_i$ if $i=D+1$. By the induction hypothesis we have $S'(E')\in\cf(V')$. 
It follows that $S(E)\in\cf(V)$, as required.

\subhead 4.5\endsubhead
Assume that $D\ge2$. Let $Sp(V)$ be the group of automorphisms of $V,(,)$. Let $\Delta$ 
be the subgroup of $Sp(V)$
generated by $R,S$ (a dihedral group of order $2(D+1)$). From 4.2(a), 4.4(a) we see that the $\Delta$-action
on $V$ induces a $\Delta$-action on $[V]$ which keeps stable the basis $\cf(V)$.

\subhead 4.6\endsubhead
We now restate the definition of $\cf(V)$ in 3.2 in more invariant terms. (In this definition the dihedral
symmetry in 4.5 is obvious.)

When $D\ge2$, we consider a connected graph with $D+1$ vertices and $D+1$ edges such that any vertex touches exactly two edges.
(this a graph of affine type $A_D$). Let $\G$ be the set of vertices and let $\L$ be the set of edges.
We assume that we are given an imbedding $\G\sub V$ such that for $\g_1\ne\g_2$ in $\G$ we have $(\g_1,\g_2)=1$ 
if $\g_1,\g_2$ are joined by an edge and $(\g_1,\g_2)=0$ if $\g_1,\g_2$ are not joined by an edge.
We then say that $(\G,\L)$ is an un-numbered circular basis (or u.c.b.) of $V$. Note that a u.c.b. exists;
in particular the circular basis $\{e_i;i\in[1,D+1]\}$ in 1.1 can be viewed as a u.c.b. in which 
$\G=\{e_i;i\in[1,D+1]\}$ and $e_i,e_j$ are joined whenever $i-j=\pm1\mod D+1$.

When $D\ge4$ we assume that $V'$ in 1.1 has a given u.c.b. with set of vertices $\G'$ and set of  edges $\L'$.
When $D\ge4$ for any $\g'\in\G',\g\in\G$ there is a unique linear map $\ti\t=\ti\t_{\g',\g}:V'@>>>V$
compatible with the symplectic forms
 and such that, setting $[\g]=\{\g\}\sqc\{\g_1\in\G;(\g_1,\g)=1\}\sub\G$, we have
$\ti\t(\g')=\sum_{\ti\g\in[\g]}\ti\g$ and $\ti\t$ restricts to a bijection $\G'-\{\g'\}@>\si>>\G-[\g]$.
This map is injective. 

We now define a collection $\cf''(V)$ of subspaces of $V$ by induction on $D$.
If $D=0$, $\cf''(V)$ consists of the subspace $\{0\}$. If $D=2$, $\cf''(V)$ consists of the subspaces of
$V$ of dimension $0$ or $1$. If $D\ge4$, a subspace $E$ of $V$ is in $\cf''(V)$ if either $E=0$ or if
there exists $g'\in\G',\g\in\G$ and $E'\in\cf''(V')$ such that $E=\ti\t_{\g',\g}(E')\op\FF_2\g$. We show:

(a) {\it If $D\ge2$ and the u.c.b. of $V$ is numbered as in 1.1 so that $\cf(V)$ is defined, we have 
$\cf''(V)=\cf(V)$.}
\nl
We argue by induction on $D$. If $D=2$ the result is obvious. Assume now that
$D\ge4$. We can assume that the u.c.b. of $V'$ is numbered as in 1.1.
For $i\in[1,D+1]$ we have

$\t_i=\ti\t_{i-1,i}$ if $2\le i\le D$,

$\t_i=\ti\t_{D-1,1}$ if $i=1$,

$\t_i=\ti\t_{D-1,D+1}$ if $i=D+1$.
\nl
Using this and the induction hypothesis we see that $\cf(V)\sub\cf''(V)$.
If $i\in[1,D-1]$ and $j\in[1,D+1]$ then for some $s\ge0$, $\ti\t_{i,j}$ is of the form $R^s\ti\t_{i,j'}$ where
$\ti\t_{i,j'}$ is as in one of the three equalities above and $R$ is as in 4.1. Using this, together with 4.2(a)
and the induction hypothesis we see that $\cf''(V)\sub\cf(V)$. This proves (a).

\head 5. Cells in Weyl groups\endhead
\subhead 5.1\endsubhead
For any finite group $\G$, let $M(\G)$ be the set consisting of pairs $(x,\r)$ where $x\in\G$ and 
$\r$ is an irreducible representation over $\CC$ of the centralizer of $x$; these pairs are taken up 
to $\G$-conjugacy; let $\CC[M(\G)]$ be the $\CC$-vector space with basis $M(\G)$ and let 
$A_\G:\CC[M(\G)]@>>>\CC[M(\G)]$ be the ``non-abelian Fourier transform'' (as in \cite{L79}).
Let $\ZZ[M(\G)]$ be the free abelian subgroup of $\CC[M(\G)]$ with basis $M(\G)$.
 
\subhead 5.2\endsubhead
In this section we fix an irreducible Weyl group $W$ and a family $c$ of irreducible representations
of $W$ (in the sense of \cite{L79}). This is the same as fixing a two-sided cell of $W$. To $c$ 
we associate a finite group $\cg_c$ as in \cite{L79}, \cite{L84}. Let $\ti\BB_c$ be the 
``new basis'' of $\CC[M(\cg_c)]$ defined in \cite{L19}. (It is actually a $\ZZ$-basis of
$\ZZ[M(\cg_c)]$.) This basis is in canonical bijection with $M(\cg_c)$, see \cite{L19}. Let
$\widehat{(x,\r)}$ be the element of $\ti\BB_c$ corresponding to $(x,\r)\in M(\cg_c)$.
We write $F$ for the non-abelian Fourier transform $A_{\cg_c}$. We have the following result.

\proclaim{Theorem 5.3} The matrix of the non-abelian Fourier transform 
$F:\CC[M(\cg_c)]@>>>\CC[M(\cg_c)]$ with respect to the new basis $\ti\BB_c$ 
is upper triangular for a suitable order on $\ti\BB_c$.
\endproclaim
From the theorem we see that there is a well defined function $\ti\BB_c@>>>\{1,-1\}$ 
(called the {\it sign function}) whose value at $\widehat{(x,\r)}\in\ti\BB_c$ is the diagonal 
entry of the matrix of $F$ at the place indexed by $\widehat{(x,\r)}$. (We use that $F^2=1$.)

In the case where $W$ is of classical type, the theorem follows from Theorem 0.2 and its proof. 
In the remainder of this section we assume that $W$ is of exceptional type.
In this case, $\cg_c$ is a symmetric group $S_n$ in $n$ letters where $n\in[1,5]$. If $n$ is $1$
or $2$ the result is immediate. The case where $n\in[3,5]$ is considered in 5.4-5.6.
We shall use the notation of \cite{L84} for the elements of $M(\cg_c)$.
Let $\th,i,\z$ be a fixed primitive root of $1$ (in $\CC$) of order $3,4,5$ respectively.

\subhead 5.4\endsubhead
In this subsection we assume that $\cg_c=S_3$.
We partition the new basis $\ti\BB_c$ in three pieces (1)-(3) as follows:
$$\widehat{(1,1)}\tag 1$$
$$\widehat{(1,r)}\tag 2$$
$$\widehat{(1,\e)},\widehat{(g_2,1)},\widehat{(g_2,\e)},
\widehat{(g_3,1)},\widehat{(g_3,\th)},\widehat{(g_3,\th^2)}.\tag 3$$
Then 

(a) {\it $F$ applied to an element in the $n$-th piece is $\pm$ that element plus a $\QQ$-linear
combination of elements in $m$-th pieces with $m>n$.}
\nl
 We have
$$\align&F(\widehat{(1,r)})\\&=F((1,1)+(1,r))=(1,1)/2+(1,r)+(1,\e)/2+(g_2,1)/2+(g_2,\e)/2\\&=
((g_2,\e)/2+(1,r)/2+(1,1)/2)+((1,\e)/2+(1,r)+(1,1)/2)+((g_2,1)/2\\&+(1,r)/2+(1,1)/2)
-((1,r)+(1,1))\\&
=\widehat{(g_2,\e)}/2+\widehat{(1,\e)}/2+\widehat{(g_2,1)}/2-\widehat{(1,r)}.\endalign$$
The formula for $F(\widehat{(1,1)})$ is as follows. If $W$ is of type $G_2$ then
$$\align&F(\widehat{(1,1)})=F(1,1)\\&=
(1,1)/6+(1,r)/3+(1,\e)/6+(g_2,1)/2+(g_2,\e)/2+(g_3,1)/3\\&+(g_3,\th)/3+(g_3,\th^2)/3\\&
=((g_3,\th)/3+(g_2,1)/3+(1,1)/3)+((g_3,\th^2)/3+(g_2,1)/3+(1,1)/3)\\&+
((g_3,1)/3+(g_2,1)/3+(1,1)/3)+((g_2,\e)/2+(1,r)/2+(1,1)/2)+((1,\e)/6\\&+(1,r)/3+(1,1)/6)
-((g_2,1)/2+(1,r)/2+(1,1)/2)-(1,1)\\&
=\widehat{(g_3,\th)}/3+\widehat{(g_3,\th^2)}/3+\widehat{(g_3,1)}/3+\widehat{(g_2,\e)}/2+
\widehat{(1,\e)}/6-\widehat{(g_2,1)}/2-\widehat{(1,1)}.\endalign$$
If $W$ is of type $E_6,E_7$ or $E_8$ then
$$\align&F(\widehat{(1,1)})=F(1,1)\\&=
(1,1)/6+(1,r)/3+(1,\e)/6+(g_2,1)/2+(g_2,\e)/2+(g_3,1)/3\\&+(g_3,\th)/3+(g_3,\th^2)/3\\&
=((g_3,\th)/3+(g_2,\e)/3+(1,1)/3)+((g_3,\th^2)/3+(g_2,\e)/3+(1,1)/3)\\&+
((g_3,1)/3+(g_2,1)/3+(1,1)/3)-((g_2,\e)/6+(1,r)/6+(1,1)/6)+((1,\e)/6\\&+(1,r)/3+(1,1)/6)
+((g_2,1)/6+(1,r)/6+(1,1)/6)-(1,1)\\&
=\widehat{(g_3,\th)}/3+\widehat{(g_3,\th^2)}/3+\widehat{(g_3,1)}/3-\widehat{(g_2,\e)}/6+
\widehat{(1,\e)}/6+\widehat{(g_2,1)}/6-\widehat{(1,1)}.\endalign$$
We see that the matrix of $F$ in the new basis is upper triangular.
This proves 5.3 in our case. 

(b) {\it The sign function on $\ti\BB_c$ is constant on each piece; its
value on the piece (1),(2),(3) is $-1,-1,1$ respectively.}

\subhead 5.5\endsubhead
In this subsection we assume that $\cg_c=S_4$ so that $W$ is of type $F_4$.
We partition the new basis $\ti\BB_c$ in five pieces (1)-(5) as follows:
$$\widehat{(1,1)}\tag 1$$    
$$\widehat{(1,\l^1)}\tag 2 $$
$$\widehat{(1,\s)} \tag 3$$  
$$\widehat{(1,\l^2)},\widehat{(g_2,1)},\widehat{(g'_2,1)},\widehat{(g_2,\e'')},\widehat{(g_2,\e')}\tag 4$$
$$\align&\widehat{(g_3,1)}, \widehat{(g_4,1)} \widehat{(g'_2,\e'')},\widehat{(g'_2,\e')},
  \widehat{(g'_2,r)},\widehat{(g_4,-1)}, \widehat{(1,\l^3)},\widehat{(g_2,\e)},\widehat{(g'_2,\e)},
\\&\widehat{(g_3,\th)},\widehat{(g_3,\th^2)},\widehat{(g_4,i)},\widehat{(g_4,-i)}.\tag 5\endalign$$
Then

(a) {\it $F$ applied to an element in the $n$-th piece is $\pm$ that element plus a $\QQ$-linear
combination of elements in $m$-th pieces with $m>n$.}
\nl
We see that the matrix of $F$ in the new basis is upper triangular. This proves 5.3 in our case.

(b) {\it The sign function on $\ti\BB_c$ is constant on each piece; its value on the piece 
(1),(2),(3),(4),(5) is $1,-1,1,-1,1$ respectively.}

\subhead 5.6\endsubhead
In this subsection we assume that $\cg_c=S_5$ so that $W$ is of type $E_8$.
We partition the new basis $\ti\BB_c$ in eight pieces (1)-(8) as follows:
$$\widehat{(g_5,\z)}\tag1$$     
$$\widehat{(1,1)}\tag2$$         
$$\widehat{(1,\l^1)}\tag3$$   
$$\widehat{(1,\nu)}\tag4$$      
$$\widehat{(1,\nu')}\tag5$$      
$$\widehat{(1,\l^2)},\widehat{(g_2,1)},\widehat{(g_2,-1)}\tag6$$   
$$\widehat{(1,\l^3)},\widehat{(g_2,r)}, \widehat{(g_3,1)},\widehat{(g'_2,1)}, 
\widehat{(g_2,-r)}, \widehat{(g'_2,r)},\widehat{(g_3,\th)},\widehat{(g_3,\th^2)}\tag7$$
$$\align&\widehat{(g'_2,\e'')},\widehat{(g_6,1)},\widehat{(g_2,\e)},\widehat{(g_3,\e)},
 \widehat{(g_4,1)},\widehat{(g_5,1)},\widehat{(g'_2,\e')},\widehat{(g_4,-1)},\\&
 \widehat{(g_6,-1)},\widehat{(g_6,\th)},\widehat{(g_6,\th^2)},\widehat{(1,\l^4)}, \widehat{(g_2,-\e)},
 \widehat{(g_3,\e\th)},\widehat{(g_3,\e\th^2)},\widehat{(g'_2,\e)},\\&
\widehat{(g_6,-\th)},\widehat{(g_6,-\th^2)},\widehat{(g_4,i)},\widehat{(g_4,-i)},\widehat{(g_5,\z^2)},
\widehat{(g_5,\z^3)},\widehat{(g_5,\z^4)}.\tag8\endalign$$
Then 

(a) {\it $F$ applied to an element in the $n$-th piece is $\pm$ that element plus a $\RR$-linear
combination of elements in $m$-th pieces with $m>n$.}
\nl
(If $n\ge2$ we can replace $\RR$ by $\QQ$ in (a). If $n=1$ the coefficients in the linear
combination can involve the golden ratio.)
We see that the matrix of $F$ in the new basis is upper triangular. This proves 5.3 in our case.

(b) {\it The sign function on $\ti\BB_c$ is constant on each piece; its value on the piece 
(1),(2),(3),(4),(5),(6),(7),(8)  is $-1,-1,1,1,1,-1,-1,1$ respectively.}
\nl
We now give some indication of how (a) can be verified.
Let $\ch$ be the hyperplane in $\CC[M(S_5)]$ consisting of all sums 
$\sum_{(x,\r)\in M(S_5)}a_{x,\r}(x,\r)$ where $a_{x,\r}\in\CC$ satisfy the equation
$$a_{g_5,\z}+a_{g_5,\z^4}=a_{g_5,\z^2}+a_{g_5,\z^3}.$$
One can check that $F(\ch)=\ch$. Moreover one can check that $\widehat{(x,\r)}\in\ch$ for any
$(x,\r)$ in $M(S_5)$ other than $\widehat{(g_5,\z)}$.
It follows that to verify (a) we can assume that $n\ge2$.
In that case the proof of (a) is similar to that of 1.6; the role of $z_i$ in 1.6 is now played by 
the maps $\ss_{H,H'}$ in \cite{L19, 3.1}; the commutation of $z_i$ with Fourier transform (see 1.5(a)) is 
replaced by the commutation of $\ss_{H,H'}$ with the non-abelian Fourier transform (see
\cite{L19, 3.1(b),(e)}). A similar argument (except for the reduction to the case $n\ge2$ which is not 
needed in this case) applies to the proof of 5.5(a).
The proof of (b) is similar to that of 2.4; we use an induction hypothesis where $S_5$ is replaced by
$S_4,S_3\T S_2,S_3,S_2\T S_2$ or $S_2$. Using the known equality $\tr(F,\CC[M_5])=13$,
we see that the values of the sign function on the elements not covered by the induction hypothesis
(that is those in the pieces (1),(2)) have sum equal to $-2$. It follows that both these values are $-1$.
A similar argument applies to the proof of 5.5(b) (in this case the only element 
not covered by the induction hypothesis is that in piece (1)).


\widestnumber\key{ABC}
\Refs
\ref\key{L79}\by G.Lusztig\paper Unipotent representations of a finite Chevalley group of type $E_8$\jour 
Quart. J. Math.\vol30\yr1979\pages315-338\endref
\ref\key{L84}\by G.Lusztig\book Characters of reductive groups over a finite field\bookinfo 
Ann.Math.Studies 107\publ Princeton U.Press \yr1984\endref
\ref\key{L19}\by G.Lusztig\paper The Grothendieck group of unipotent representations: a new basis\lb\jour
arxiv:1907.01401\endref
\endRefs
\enddocument